\documentclass[12pt]{article}
\usepackage{amsmath}
\usepackage{amssymb}
\usepackage{dsfont}
\usepackage{cite}
\newsymbol \blackbox 1004
\newcommand{\eh}{\hfill}\newlength{\sperr}

\newenvironment{proof}{{\settowidth{\sperr}{\bf\rm
Proof}%
\par\addvspace{0.3cm}\noindent\parbox[t]{1.3\sperr}
{\bf\rm P\eh r\eh o\eh o\eh f\eh }%
}}{\nopagebreak\mbox{}
$\blackbox$\par\addvspace{0.3cm}}

\def\nn{\nonumber}

\def\a{\alpha}
\def\b{\beta}

\def\g{\gamma}

\def\de{\delta}
\def\De{\Delta}

\def\vk{\varkappa}

\def\Lam{\Lambda}
\def\s{\sigma}
\def\la{\lambda}
\def\om{\omega}

\def\th{\theta}

\def\up{\upsilon}
\def\vp{\varphi}
\def\vt{\vartheta}

\def\wh{\widehat}
\def\wt{\widetilde}
\def\ov{\overline}

\def\BC{{\mathbb C}}
\def\BD{{\mathbb D}}
\def\BR{{\mathbb R}}
\def\BN{{\mathbb N}}
\def\BT{{\mathbb T}}
\def\clp{{\mathcal P}}

\def\clb{{\mathcal B}}

\def\cld{{\mathcal D}}
\def\cle{{\mathcal E}}
\def\clf{{\mathcal F}}
\def\clg{{\mathcal G}}
\def\clh{{\mathcal H}}

\def\clk{{\mathcal K}}
\def\clm{{\mathcal M}}
\def\cln{{\mathcal N}}
\def\clu{{\mathcal U}}

\def\clr{\mathcal{R}}

\def\clw{{\mathcal W}}
\def\cld{{\mathcal D}}

\newcommand{\E}{\mathrm{e}}
\newcommand{\I}{\mathrm{i}}

\def\Res{\mathrm{Res}}

\newtheorem{Pa}{Paper}[section]
\newtheorem{Tm}[Pa]{{\bf Theorem}}
\newtheorem{La}[Pa]{{\bf Lemma}}

\newtheorem{Rk}[Pa]{{\bf Remark}}

\newtheorem{Dn}[Pa]{{\bf Definition}}
\newtheorem{Nn}[Pa]{{\bf Notation}}
\newtheorem{Pn}[Pa]{{\bf Proposition}}

\newenvironment{dedication}
        {\vspace{1ex}\begin{quotation}\begin{center}\begin{em}}    
        {\par\end{em}\end{center}\end{quotation}}

\title{Arov--Krein entropy functionals and indefinite interpolation problems}

\author{I. Roitberg and A.L.~Sakhnovich}

\date{}
\begin{document}
\maketitle

\begin{dedication} 
\end{dedication}

\begin{abstract}  We generalize the notion of the Arov-Krein entropy functional for the case
of generalized Nevanlinna functions and obtain a representation of these functionals on solutions
of indefinite interpolation problems. The case of indefinite Caratheodory problem
and application to Szeg{\H{o}} limit formula for this nonclassical case are considered
in greater detail.
\end{abstract}


\vspace{0.2em}


\section{Introduction}\label{Intro}
\setcounter{equation}{0}
{\bf 1.} In their important work \cite{ArK}, D.Z. Arov and M.G. Krein studied matrix-valued functions $\om(\la)$ (given as linear fractional transformations),
such that $\om(\la)+\om(\la)^*\geq 0$  on the unit disk $\BD=\{\la: \, |\la |<1\}$, where $\om(\la)^*$ stands for the complex conjugate transpose of $\om(\la)$.  (See also related results 
in \cite{ArK2, ALS87}.)  The class of such matrix-valued functions (matrix functions) is called Caratheodory class $C$ (or $C_0$) and $\om\in C$
admits Herglotz representation 
\begin{align} \label{I1}
\om(\la)=\I \nu+\frac{1}{2\pi}\int_{0}^{2\pi}\frac{\E^{\I\th}+\la}{\E^{\I\th}-\la}d\mu(\th) \quad (\nu=\nu^*),
\end{align}
where $\mu(\th)$ is a nondecreasing matrix function. In particular, one can easily see that the real part
\begin{align} \label{I1+}
\Re\left(\frac{\E^{\I\th}+\la}{\E^{\I\th}-\la}\right)=\frac{1-|\la |^2}{|\E^{\I \th}-\la |^2}\geq 0.
\end{align}
The Arov--Krein entropy functionals are given by the formula:
\begin{align} \label{I3}
E(\om, \wt \la)=-\frac{1}{4\pi}\int_{0}^{2\pi}|\E^{\I\th}-\wt \la |^{-2}(1-| \wt\la |^2)\ln \det\big(\mu^{\prime}(\th)\big)d \th \quad (|\wt \la |<1).
\end{align}
Further results on Arov-Krein and related Burg's entropy and further references are given in \cite{Bult, Burg, DymG, FFK, Leon}.

Here, we consider a natural generalization of  the Arov--Krein entropy functional for the case of the solutions of indefinite
interpolation problems belonging to the Krein--Langer (generalized Caratheodory) class $C_{\vk}$, where $C_0=C$. 
More precisely, we take solutions of indefinite
interpolation problems belonging to the generalized Nevanlinna class $N_{\vk}$, transform them into $C_{\vk}$
and consider the entropy.
See Definition \ref{DnN} of the classes $N_{\vk}$ and $C_{\vk}$ and a more detailed discussion
below. The problem is closely related to the Szeg{\H{o}} limit formula in the indefinite case \cite{ALSind, ALSJFA} (see also
interesting recent papers \cite{DerSi, DerSi2}).

{\bf 2.} The solutions of the so called {\it indefinite interpolation problems} are often described in terms of linear fractional transformations.
In this paper, we deal with the case when the matrix functions, which are obtained using these linear  fractional transformations, are defined
in the open upper  half-plane $\BC_+$.  Moreover, these functions belong to the classes $N_{\vk}$, where $\vk \geq 1$ (see the
definition of $N_{\vk}$ below). 
\begin{Dn}\label{DnN}
The generalized Nevanlinna class $N_{\vk}$ is the set of  meromorphic $p \times p$  matrix functions $\vp(z)$  $(z\in \BC_+)$ such that the kernel
$\big(\vp(z) - \vp(\zeta)^*\big)\big/ (z-\ov{\zeta})$ has $\vk$ negative squares. That is, for any $n \in \BN$ and any set $z_1, \, z_2,
 \, \ldots , \, z_n \, \in \, \BC_+$ the matrix $\big\{\big(\vp(z_i) - \vp(z_k)^*\big)\big/ (z_i-\ov{z_k})\big\}_{i,k=1}^n$ has at most $\vk$ negative eigenvalues
 and for at least one choice $n, \, z_1, \, z_2,  \, \ldots , \, z_n$ it has exactly $\vk$ negative eigenvalues.
 
 The generalized Caratheodory class $C_{\vk}$ is the set of  meromorphic $p \times p$  matrix functions $\om(\la)$  $(|\la |<1))$ such that the kernel
$\big(\om(\la) + \om(\zeta)^*\big)\big/ (1-\la \ov{\zeta})$ has $\vk$ negative squares. That is, for any $n \in \BN$ and any set $\la_1, \, \la_2,
 \, \ldots , \, \la_n \, \in \, \BD$ the matrix $\big\{\big(\om(\la_i) + \om(\zeta_k)^*\big)\big/ (1-\la_i \ov{\zeta_k}\big\}_{i,k=1}^n$ has at most $\vk$ negative eigenvalues
 and for at least one choice $n, \, \la_1, \, \la_2,  \, \ldots , \, \la_n$ it has exactly $\vk$ negative eigenvalues.
\end{Dn}
Here $\BN$ denotes the set of natural numbers and $\ov{z}$ is the complex conjugate of $z$.
The classes $N_{\vk}$ and $C_{\vk}$ have been studied in a series of seminal Krein--Langer papers, for instance \cite{KrL, KrL2} (see also the articles
\cite{DaLang, RoS0} and references
therein). Clearly, $C_0=C$.
\begin{Rk} Sometimes, we use the notations $N_{\vk}^{p}$ and $C_{\vk}^{p}$ instead of $N_{\vk}$ and $C_{\vk}$, respectively,
 to stress the order $p$ of the matrix functions in the classes.
\end{Rk}

We will use the ``operator identities" formalism of solving  indefinite interpolation problems from \cite{ALSind} (scalar function case)
and \cite{RoS1, RoS2, ALSJFA} (matrix function case), which is based on the approaches to sign-definite interpolation problems
described in \cite{Pot1, SaL2}. This formalism as well as the representation of the generalized Nevanlinna  functions is
discussed in Section ``Preliminaries" (Section \ref{Prel}). In Section \ref{Ent}, we introduce entropy functional for
the generalized Nevanlinna  functions and obtain a representation of this functional on the solutions of  indefinite
interpolation problems. In Section \ref{Car}, we consider in detail the case of block Toeplitz matrices and indefinite
Caratheodory problem. Section \ref{lim} is dedicated to the application of the results from Section \ref{Car}  in
the proof of indefinite Szeg{\H{o}} limit formula.

 In the paper, $\BN$ denotes the set of natural numbers, $\BN_0$ denotes the set of nonnegative integers,
 $\BR$ denotes the real axis,  $\BD$ is the unit disk (i.e., $\BD=\{\la: \, |\la |<1\}$),
 $\BT$ stands for the unit circle (i.e., $\BT=\{\la: \, |\la |=1\}$), $\BC$ stands for the complex plane, and
$\BC_+$ ($\BC_-$) stands for the open upper (lower) half-plane. For $\De \subset \BR$, $\chi_{\De}$ is the characteristic function
of $\De$, that is, $\chi_{\De}(t)=1$ for $t\in \De$ and $\chi_{\De}(t)=0$ for $t\not\in \De$.
The spectrum of a bounded operator or a square matrix $A$ is denoted by $\s(A)$ and $I_p$ stands for the $p\times p$ identity matrix.

The set of linear bounded operators acting from  the Hilbert space $\clg$ into Hilbert space $\clh$ is denoted by $\clb(\clg,\clh)$ and the notation
$\clb(\clh, \clh)$ we simplify into $\clb(\clh)$.
The set of Hermitian operators (or matrices) $S$
from $\clb(\clh)$, such that the spectrum of $S$ contains precisely  $\vk$ (counting multiplicities) negative eigenvalues ($\vk<\infty$), is denoted by
$\clp_{\vk}$. We use also the notation $\vk_S$ for the mentioned above $\vk$.

\section{Preliminaries }\label{Prel}
\setcounter{equation}{0}
\subsection{Representation of $\vp(z)\in N_{\vk}$}
The matrix function $\vp(z)\in N_{\vk}^p$ admits representation \cite{DaLang} (see also \cite{RoS0}):
\begin{align}  \label{1}  
\vp(z)=&\int_{-\infty}^{\infty}\left(\frac{1}{t-z}+\sum_{i=0}^rK_i(t,z)\right)d\tau(t)+R_0(z)
\\ \nn &
- \sum_{i=1}^rR_i\left(\frac{1}{z-\a_i}\right)+\sum_{j=1}^m\left(M_j\left(\frac{1}{z-\b_j}\right)+M_j\left(\frac{1}{ \ov{z}-\b_j}\right)^*\right),
\end{align} 
where $\a_1<\a_2<\ldots <\a_r$ are real numbers, $\b_1, \ldots , \b_m \in \BC_+$ are distinct numbers;\\
(i)  the real line $\BR$ is a union of the sets $\De_0, \, \De_1, \ldots, \De_r$, such that $\De_1, \ldots, \De_r$ are bounded open intervals
having disjoint closures, $\De_0$ is there complement, and $\a_i\in \De_i$ $(1 \leq i \leq r)$;
\begin{align}  \label{2}  &
K_i(t,z)=\sum_{k=1}^{2\rho_i}\frac{(t-\a_i)^{k-1}}{(z-\a_i)^{k}}\chi_{\De_i}(t) \quad  {\mathrm{for}} \quad 1 \leq i \leq r
\end{align} 
(recall the definition of the characteristic function $\chi_{\De}$ from the introduction),
\begin{align}  \label{2'}  &
K_0(t,z)=-\left(t\frac{(1+z^2)^{\rho_0}}{(1+t^2)^{\rho_0+1}}+(t+z)\sum_{k=1}^{\rho_0}\frac{(1+z^2)^{k-1}}{(1+t^2)^k}\right)\chi_{\De_0}(t),
\\ \nn & \rho_0\in \BN_0, \quad \rho_1, \, \ldots , \, \rho_r \, \in \, \BN;
\end{align} 
(ii) $\tau(t)$ is a nondecreasing (on each of the intervals $(-\infty, \a_0)$,
$(\a_i,\a_{i+1})$, where $0\leq i<r$,  and $(\a_r, \infty)$) $p \times p$ matrix function such that  the following integral converges:
\begin{align}  \label{2+}  &
\int_{-\infty}^{\infty}\frac{(t-\a_1)^{2\rho_1}\ldots (t-\a_r)^{2\rho_r}}{(1+t^2)^{\rho_1+\ldots +\rho_r}}\frac{d \, \tau(t)}{(1+t^2)^{\rho_0+1}}<\infty;
\end{align}
(iii) for each $0 \leq i \leq r$, the $p\times p$ matrix function $R_i(z)$ is a matrix polynomial of degree at most $2\rho_i+1$
having self-adjoint matrix coefficients, such that the coefficient $C_i$ in the term $C_i z^{2\rho_i+1}$ of maximal degree
in $R_i$ is  nonnegative ($C_i \geq 0$), and the equalities $R_1(0)= \cdots = R_r(0)=0$ hold; \\
(iv) for each $1 \leq j \leq m$, the $p\times p$ matrix function $M_j(z)$ is a matrix polynomial $\not\equiv 0$ such that  $M_j(0)=0$.
\begin{Rk}\label{RkSum} It is easy to see that
\begin{align}  \label{2=}  &
\frac{1}{t-z}+K_i(t,z)=\frac{1}{t-z}\left(\frac{t-\a_i}{z-\a_i}\right)^{2\rho_i} \quad  {\mathrm{for}} \quad t\in \De_i \quad (1 \leq i \leq r),
\\  \label{2/}  &
\frac{1}{t-z}+K_0(t,z)=\frac{1+tz}{t-z}\frac{(1+z^2)^{\rho_0}}{(1+t^2)^{\rho_0+1}} \quad  {\mathrm{for}} \quad t\in \De_0.
\end{align} 
\end{Rk}
\begin{Nn} The degree of $M_j$ is denoted by $\zeta_j$. \end{Nn}
\begin{Rk}\label{Rkp1}
For $p=1$,   without loss of generality $($see \cite[Theorem 3.1]{KrL}$)$ we  require that
\begin{align}  \label{2!}  &
\sum_{i=0}^r \rho_i +\sum_{j=1}^m \zeta_j=\vk.
\end{align}
\end{Rk}  
\subsection{Indefinite interpolation problem}
Let the operator (or matrix) $S\in \clp_{\vk}$ be given.
In the following text we assume that for some $A\in \clb(\clh)$ and $p\in \BN$ the operator identity
\begin{align}  \label{3}  &
AS-SA^*=\I (\Phi_1\Phi_2^*+\Phi_2\Phi_1^*) \qquad \big(\Phi_1\in \clb(\BC^p, \clh), \quad \Phi_2 \in \clb(\BC^p, \clh)) 
\end{align} 
holds. Fixed operators $A$ and $\Phi_2$ determine a class of so called structured operators satisfying
\eqref{3} (e.g., Toeplitz or Loewner matrices, operators with difference kernels and so on). 
We assume also that $S$ is invertible and
\begin{align}  \label{4}  &
 S^{-1}\in \clb(\clh).
\end{align} 
{\bf 1.} We  introduce the operators $S_{\vp}$ and $\Phi_{\vp}$ in terms of the representation \eqref{1} of $\vp\in N_{\vk}^p$.
For this purpose we need some preparations.
Let $A$ and $\Phi_2$ be fixed and  the matrix function $F_j$ and $\g_j \in \BC$ $(j>0)$ be given. Then we set
\begin{align}  \label{11}  &
\clf_j=\Res_{z=\g_j}(I-zA)^{-1}\Phi_2F_j(z)\Phi_2^*(I-zA^*)^{-1} \quad (j>0),
\\ \label{12}  &
\wh \clf_j=\Res_{z=\g_j}(I-zA)^{-1}A\Phi_2F_j(z) \quad (j>0);
\\ \label{13}  &
\clf_0=\Res_{z=0}(A-zI)^{-1}\Phi_2F_0(1/z)\Phi_2^*(A^*-zI)^{-1},
\\ \label{14}  &
\wh \clf_0=-\Res_{z=0}\Big((A-zI)^{-1}A\Phi_2F_0(1/z)\big/ z\Big) .
\end{align} 
\begin{Dn}\label{DnKLD}
The Krein--Langer data corresponding to the representation \eqref{1} is the set
\begin{align}  \nn
\cld=&\big\{\tau(t); \,\, \a_1,\ldots \a_r; \,\, \b_1, \ldots, \b_m; \,\, \rho_0, \rho_1, \ldots \rho_r; \,\, \De_0, \De_1, \ldots , \De_r; 
\\ \label{KLD}  & \quad
R_0, R_1, \ldots, R_r; \,\, M_1, \ldots, M_m\big\}.
\end{align} 
\end{Dn}
Clearly the data $\cld(\vp)$ corresponding to $\vp\in N_{\vk}$ is not defined uniquely although the arbitrariness  is not so great
(see Remark \ref{RkSvp}). In particular,  the points $\b_1, \ldots, \b_m$ and the functions $M_1, \ldots, M_m$
are fixed. The following operators and operator functions are generated by the data $D(\vp)$
using \eqref{11}--\eqref{14}.
\begin{Dn}\label{DnF} When $F_j(z)=R_j\big(1\big/ (z-\a_j)\big)$ and $\g_j=\a_j$, we denote the corresponding $\clf_j$ and $\wh \clf_j$  by $\clr_j$ and $\wh \clr_j$,
respectively;
when $F_j(z)=M_j\big(1\big/ (z-\b_j)\big)$ and $\g_j=\b_j$, we denote the corresponding $\clf_j$ and $\wh \clf_j$  by $\clm_j$ and $\wh \clm_{1j}$,
respectively, and when $F_j(z)=M_j\big(1\big/ (\ov{z}-\b_j)\big)^*$ and $\g_j=\ov{\b_j}$ we denote  the corresponding $\wh \clf_j$   by $\wh \clm_{2j}$.
When $F_j(z)=K_j(t,z)$ and $\g_j=\a_j$, we denote the corresponding  $\wh \clf_j$ by $\wh \clk_j(t,z)$.
When $F_0(z)=R_0(z)$, we denote the corresponding $\clf_0$ and $\wh \clf_0$  by $\clr_0$ and $\wh \clr_{0}$,
respectively. Finally, setting $F_0(z)=K_0(t,z)$, where $t$ is an additional parameter, we put
$$\wh \clk_0(t):=-\clf_0(t).$$
The matrix functions $\tau_j(t)$ are introduced by the similar to \eqref{11} formulas$:$
\begin{align}&\label{15}
d\tau_j(t)=\Res_{z=\a_j}\Big(K_j(t,z)(I-zA)^{-1}\Phi_2\big(d\tau(t)\big)\Phi_2^*(I-zA^*)^{-1} \Big)
\end{align}
for $j>0$, and
\begin{align}
 &\label{16}
d\tau_0(t)=-\Res_{z=0}\Big(K_0\big(t,1/z\big)(A-zI)^{-1}\Phi_2\big(d\tau(t)\big)\Phi_2^*(A^*-zI)^{-1} \Big).
\end{align}
\end{Dn}
Let $A$, $\Phi_2$ and the representation \eqref{1} of $\vp(z)$ be given. Then, taking into account Definition \ref{DnF},
we introduce the operators
\begin{align}  \nn
S_{\vp}=& \int_{-\infty}^{\infty}\left((I-tA)^{-1}\Phi_2 \big(d\tau(t)\big)\Phi_2^*(I-tA^*)^{-1}-\sum_{i=0}^rd\tau_j(t)\right)
\\ & \label{20}  
+\sum_{i=0}^r \clr_j-\sum_{j=1}^m\big(\clm_j+\clm_j^*\big);
\\ \nn 
\I \Phi_{\vp}=&\int_{-\infty}^{\infty}\left((I-tA)^{-1}A\Phi_2 -\sum_{i=0}^r\wh \clk_i(t)\right)d\tau(t)
\\ & \label{21}  
+\sum_{i=0}^r\wh \clr_i-\sum_{j=1}^m\big(\wh \clm_{1j}+\wh \clm_{2j}\big).
 \end{align} 

\begin{Rk} \label{RkSvp} According to \cite[Section 3]{RoS1} the operators $S_{\vp}$ and $\Phi_{\vp}$ are well-defined
and  the integrals in \eqref{20} and \eqref{21} weakly converge
under conditions \eqref{9} and 
\begin{align}  \label{10}  &
1\big/\b_j \not\in \s(A), \quad 1\big/ \, \ov{\b_j} \not\in \s(A).
\end{align} 
Moreover, Theorems 3.4 and 3.5 in \cite{RoS1} state that 
$S_{\vp}$ and $\Phi_{\vp}$ do not depend on the choice of the domains $\De_0, \ldots \De_r$, that
$S_{\vp}=S_{\vp}^*\in \clp_{\kappa}$, where $\kappa <\infty$, and the operator identity
\begin{align}  \label{3'}  &
AS_{\vp}-S_{\vp}A^*=\I (\Phi_{\vp}\Phi_2^*+\Phi_2\Phi_{\vp}^*) 
\end{align} 
holds.
\end{Rk}
{\bf 2.} The transfer matrix function in Lev Sakhnovich form \cite{SaSaR, SaL1, SaL2} is introduced by the equality
\begin{align}  \label{5}  &
w_A(z)=I_{2p}+\I zJ\Pi^* S^{-1}(I-zA)^{-1}\Pi, \quad \Pi:=\begin{bmatrix}\Phi_1 & \Phi_2\end{bmatrix}, \quad J:=\begin{bmatrix}0 & I_p \\ I_p & 0\end{bmatrix},
\end{align} 
where $A$, $S$ and  $\Pi$ satisfy \eqref{3}, $I$ is the identity operator and $I_p$ is the $p \times p$ identity matrix. We will also need the matrix function
\begin{align}  \label{6}  &
\clu(z):=w_A(\ov{z})^*=I_{2p}-\I z\Pi^* (I-zA^*)^{-1}S^{-1}\Pi J= \begin{bmatrix}a(z) & b(z) \\ c(z) & d(z) \end{bmatrix},
\end{align} 
where $a,b,c$ and $d$ are $p \times p$ blocks of  $\,  \clu$. By virtue of the properties of the transfer matrix functions 
(see, e.g. \cite[(1.84)]{SaSaR}) we have a useful equality 
\begin{align} \label{En2-}
\clu(\ov{z})^*J\clu(z)\equiv J.
\end{align}
The linear fractional transformations, which we
are interested in, are given by the formula
\begin{align}  \label{7}  &
\vp(z)=\vp(z,P,Q)=\I(a(z)P(z)+b(z)Q(z))(c(z)P(z)+d(z)Q(z))^{-1},  \\
\label{7'}  &
 \det\big(c(z)P(z)+d(z)Q(z)\big)\not\equiv 0,
 \end{align} 
where  the $p \times p$ matrix functions $P$ and $Q$ are meromorphic in $\BC_+$  and satisfy the inequalities
\begin{align}  \label{8}  &
P(z)^*Q(z)+Q(z)^*P(z)\geq 0, \quad P(z)^*P(z)+Q(z)^*Q(z)> 0,
 \end{align} 
except at isolated points. One says that such pairs $\{P,Q\}$   are {\it nonsingular  pairs with the property-$J$}.
Given a so called {\it frame} $\clu$, let us introduce $\cln(\clu)$.
\begin{Nn} \label{LFT} The set of matrix  functions $\vp(z)$ of the form \eqref{7},
where $\{P,Q\}$   are  nonsingular  pairs with the property-$J$
$($which satisfy \eqref{7'}$)$, is denoted by $\cln(\clu)$. \end{Nn}

In this paper, we consider the simpler interpolation cases, where the following conditions are valid:
\begin{align}  \label{9}  &
\s(A)\cap \s(A^*)= \emptyset ; \qquad \s(A) \quad {\mathrm{is \,\, a \,\, finite \,\, set}}.
 \end{align} 

{\bf 3.}    The function $\vp(z)$ given in $\BC_+$ is determined (further in the text) in $\BC_-$ by the 
relation 
\begin{align}  \label{vp}  &
\vp(z):=\vp(\ov{z})^* \quad {\mathrm{for}} \quad z\in \BC_-.
 \end{align} 
 Now, we introduce  ${\bf B}_{\vp}(z)$ in $\BC_+\cup \BC_-$ by the equality
\begin{align}  \label{22}  &
{\bf B}_{\vp}(z):=(I-zA)^{-1}\big(\Phi_1-\I \Phi_2 \vp(z)\big).
\end{align} 
\begin{Nn} \label{NnE} We denote by $\cle$  the class of matrix functions
$\vp(z)$ which are analytic in $\BC_+\cup \BC_-$ $($excluding, possibly, isolated points$)$
and may have only removable singularities in the points $z$ such that 
$1/z \in \s(A)$.
\end{Nn}
The next interpolation theorem directly follows from \cite[Theorems 4.4, 5.1]{RoS1}.
\begin{Tm}  \label{TmInt}Let the operators  $A$, $S$, $\Phi_1$ and $\Phi_2$
satisfy the operator identity \eqref{3}.  Assume that $S\in \clp_{\vk_S}$, that
conditions \eqref{4} and \eqref{9} are fulfilled and that $\clu$ is given by \eqref{6}.
Then the following two statements are valid. \\
$(i)$ If $\vp\in \cln(\clu)$ $($in $\BC_+)$, $\vp \in \cle$ $($see Notation \ref{NnE}$)$   and the function
${\bf B}_{\vp}(z)$
is analytic at every $z$ such that $1/z \in \s(A)$, then $\vp\in N_{\vk}$ $(\vk \leq \vk_S)$, \eqref{10} is valid and the equalities
\begin{align}  \label{23}  &
S=S_{\vp}, \quad \Phi_1=\Phi_{\vp},
\end{align}
where $S_{\vp}$ and $\Phi_{\vp}$ are given by \eqref{20} and \eqref{21}, respectively, are satisfied. \\
$(ii)$ Conversely, if some $p \times p$ matrix function $\vp(z)$ belongs $N_{\vk}$ and relations \eqref{10} and \eqref{23} are fulfilled,
then $\vp(z)\in \cln(\clu)$, $\vp(z)\in \cle$ and ${\bf B}_{\vp}(z)$ is analytic at every $z$ such that $1/z \in \s(A)$.
\end{Tm}
Clearly, one can consider the linear fractional transformation \eqref{7} in both half-planes $\BC_+$ and $\BC_-$ as is done
in \cite{RoS1, ALSind} and is explained for the case $p=1$ below.
\begin{Dn} \label{DnSIP} Let the conditions of Theorem \ref{TmInt} hold.
Then, the matrix functions $\vp(z)$ considered in  Theorem \ref{TmInt}
are called the solutions of the interpolation problem \eqref{23}.
\end{Dn}
In the scalar case $p=1$,  the pairs $\{P,Q\}$ look somewhat simpler (see the remark below and \cite{ALSind}).
\begin{Rk}
Assuming $p=1$ in this remark, we note that without loss of generality one may consider linear fractional transformations
\eqref{7} generated by  the pairs $P(z)=\psi(z)\in N$ and $Q(z)\equiv \I$
completed with the pair $P(z)\equiv 1$ and $Q(z)\equiv 0$. In other words, we have
\begin{align}  \label{24}  
\cln(\clu)=\{\vp(z): \, &\vp(z)=\I(a(z)\psi(z)+\I b(z))(c(z)\psi(z)+\I d(z))^{-1}, 
\\ \nn &
 \psi(z)\in N, \quad c(z)\psi(z)+\I d(z)\not\equiv 0\}\cup \cln_{\infty},
\end{align}
where
\begin{align}  \label{25}  &
\cln_{\infty}=\{\I a(z)\big/c(z)\} \quad {\mathrm{if}} \quad c(z)\not\equiv 0; \quad \cln_{\infty}=\emptyset \quad {\mathrm{if}} \quad c(z) \equiv 0.
\end{align} 
\end{Rk}
We set $\psi(\ov{z})=\ov{\psi(z)}$ . Formula \eqref{En2-} may be rewritten in the form
\begin{align}  \label{26}  &
\clu(\ov{z})=(J\clu(z)^*J)^{-1}=J\big(\clu(z)^{-1}\big)^*J.
\end{align} 
When $p=1$, formula \eqref{26} takes the form
\begin{align}  \label{27}  &
\clu(\ov{z})=\big(1\big/ \, \ov{\det(\clu(z))} \, \big)\begin{bmatrix} \ov{a(z)} &-\ov{b(z)}\\ -\ov{c(z)} & \ov{d(z)}\end{bmatrix}.
\end{align} 
Recall that 
\begin{align}  \label{28}  &
\vp(z)=\ov{\vp\big(\ov{z}\big)}, \quad \psi(z)=\ov{\psi\big(\ov{z}\big)} \quad {\mathrm{for}} \quad z\in \BC_-.
\end{align} 
In view of \eqref{27} and \eqref{28}, the functions given in \eqref{24} for $z\in \BC_+$ have the same form in $\BC_-$, that is
\begin{align}  \label{29}  &
\vp(z)=\I(a(z)\psi(z)+\I b(z))(c(z)\psi(z)+\I d(z))^{-1} \\
\nn &
 {\mathrm{or}} \quad \vp(z)=\I a(z)\big/c(z) \quad {\mathrm{for}} \quad z\in \BC_+\cup \BC_-.
\end{align} 
Instead of the condition $\vp(z)\in \cle$ in Theorem \ref{TmInt}, we will use the condition
$c(z)\psi(z)+\I d(z)\not=0$ when $1/z\in \s(A)$ in the next theorem.

\begin{Tm} \label{TmIntp1} Let $p=1$ and let the operators  $A$, $S$, $\Phi_1$ and $\Phi_2$
satisfy the operator identity \eqref{3}.  Assume that $S\in \clp_{\vk_S}$ and that
conditions \eqref{4} and \eqref{9} are fulfilled.
Then the following two statements are valid.\\
$(i)$ If $\vp\in \cln(\clu)$  and $c(z)\psi(z)+\I d(z)\not=0$ $(c(z)\not=0$ when $\vp=\I a(z)\big/c(z))$ for all $z$ such that $1/z\in\s(A)$, then
$\vp\in N_{\vk_S}$  and
the equalities \eqref{10} and  \eqref{23} are satisfied.\\
$(ii)$ Conversely, if some function $\vp(z)$ belongs $N_{\vk_S}$ and relations \eqref{10} and \eqref{23} are fulfilled,
then $\vp(z)\in \cln(\clu)$.
\end{Tm}

\section{Entropy}\label{Ent}
\setcounter{equation}{0}
There is a simple one to one mapping between the matrix functions  $\vp(z)\in N_{\vk}^p$ and matrix functions $\om(\la) \in C_{\vk}^p$
(see, e.g., \cite{KrL} or \cite[p. 344]{RoS2}):
\begin{align} \label{NC}
\om (\la)=-\I \vp(z), \quad {\mathrm{where}} \quad \la=\frac{z-\up}{z-\ov{\up}} \quad (\up \in \BC_+).
\end{align}
Correspondingly, the representation \eqref{1} of $\vp(z)$ is equivalent to the following representation of $\om(\la)$:
\begin{align}  \label{1'}  
\om(\la)=&\frac{1}{2\pi}\int_{(0,2\pi)\setminus\{a_1, \ldots,a_r\}}\left(\frac{\E^{\I\th}+\la}{\E^{\I\th}-\la}+\sum_{i=0}^{r}L(a_i,\rho_i, \th,\la)\chi_{\wt \De_i}\right)d\mu(\th)+T(\la),
\end{align} 
where $\mu(\th)$ is a nondecreasing on the intervals $(0,a_1)$, $(a_i,a_{i+1})$ for $0<i<r$ and $(a_r,2\pi)$ matrix function,  
and the following integral converges:
\begin{align} \label{NC0-}
\int_{(0,2\pi)\setminus\{a_1, \ldots,a_r\}}\th^{2\rho_0}(\th-a_1)^{2\rho_1}\cdots (\th-a_r)^{2\rho_r}(2\pi -\th)^{2\rho_0}d\mu(\th)<\infty;
\end{align}
$\E^{\I a_0}=1$ (i.e., $a_0=0$ or $a_0=2\pi$);
 $\E^{\I a_i}= \frac{\a_i-c}{\a_i-\ov{c}}$ $(0<a_i<2\pi)$ for $i>0$,  where the set ${\a_i}$ is taken from the representation
\eqref{1} of $\vp(z)$; \\ $\chi_{\wt \De_i}$ is a characteristic function and for $i>0$  the set $\wt \De_i$ is an open interval on $[0,2\pi]$ containing $\a_i$;
 $\wt \De_0=[0,2\pi]\setminus \bigcup_{i=1}^r\wt \De_i$;
\begin{align} \label{NC0}
T(\la)=&-\I R_0(z)
+\I \sum_{i=1}^rR_i\left(\frac{1}{z-\a_i}\right)
\\ \nn &
-\I \sum_{j=1}^m\left(M_j\left(\frac{1}{z-\b_j}\right)+M_j\left(\frac{1}{ \ov{z}-\b_j}\right)^*\right);
\end{align}
for $ \E^{\I a}\not=1$ we have
\begin{align}& \label{NC1}
L(a, \rho, \th, \la)=\frac{2\E^{\I \th}(\E^{\I a}-1)}{(\E^{\I \th}-1)^2}\sum_{k=0}^{2\rho-1}\left(\frac{\E^{\I \th}-\E^{\I a}}{\E^{\I \th}-1}\right)^k
\left(\frac{\la-1}{\la-\E^{\I a}}\right)^{k+1}-\frac{\E^{\I \th}+1}{\E^{\I \th}-1}, 
\end{align}
and for $ \E^{\I a}=1$ we have
\begin{align}& \label{NC2}
L(a, \rho, \th, \la)=\frac{(\E^{\I \th}+\la)(\la \E^{\I \th}-1)}{\E^{\I \th}(1-\la)^2}\sum_{k=0}^{\rho-1}\frac{\la^k(\E^{\I \th}-1)^{2k}}{\E^{\I k\th}(\la-1)^{2k}}.
\end{align}
According to \cite[formula (2.10)]{RoS2} we have
\begin{align}& \label{NC3}
\frac{\E^{\I \th}+\la}{\E^{\I \th}-\la}+L(a, \rho, \th, \la)=O(\th-a)^{2\rho} \quad {\mathrm{for}} \quad \th \to a \quad (\la \in \BD),
\end{align}
where $a$ and $\rho$ take the values $a_i$ and $\rho_i$ $(1\leq i \leq r)$, respectively, as well as $a=0$, $\rho=\rho_0$
and $a=2\pi$, $\rho=\rho_0$.  We have also  (see \cite[formula (2.9)]{RoS2}) a useful relation
\begin{align}& \label{NC4}
\ov{L(a, \rho, \th, \la)}=-L(a, \rho, \th, 1\big/\ov{\la}).
\end{align}
\begin{Dn}\label{ENk} The entropy functionals on the matrix function $\vp\in N_{\vk}$
are given by the  formula
\begin{align} \label{E0}
E(\vp,\wt \la)=-\frac{1}{4\pi}\int_{0}^{2\pi}|\E^{\I\th}-\wt \la |^{-2}(1-|\wt \la |^2)\ln \det\big(\mu^{\prime}(\th)\big)d \th \quad (|\wt \la|<1),
\end{align}
where $\mu$ is the matrix function from the Krein--Langer representation \eqref{1'} of $\om(\la)$
and $\om$ and $\la$ are constructed from $\vp$ and $z$, respectively, using \eqref{NC}. 
\end{Dn}
Note that the right-hand sides of \eqref{I3} and \eqref{E0} formally coincide although the requirements on $\mu$ 
in \eqref{I3} and \eqref{E0} differ.

Consider representation \eqref{1'} in greater detail. Fix some arbitrary closed interval $[\ell_1, \ell_2]$  belonging to one of the intervals $(0, \a_1)$,
$(\a_i,\a_{i+1})$, where $0< i<r$,  or $(\a_r, 2\pi)$, and note that  the matrix function
$\mu(\th)$ is bounded on $[\ell_1, \ell_2]$. Hence, we write down $\om(\la)$ as a sum of two functions generated by the representation \eqref{1'}:
\begin{align}  \label{E1}  &
\om(\la)=\om_0(\la)+\wt \om(\la), \quad \wt \om(\la):=\frac{1}{2\pi} \int_{\ell_1}^{\ell_2}\frac{\E^{\I\th}+\la}{\E^{\I\th}-\la}d\mu(\th).
 \end{align} 
It is easy to see that  $\wt \om\in C$. Taking into account Smirnov's and Nevanlinna's theorems, we derive that the entries
of $\wt \om(\la)$ belong to the Hardy class $H_{\de}$ for each $1>\de>0$, and the non-tangential limits  $\lim_{\la \to \exp\{\I t\}}  \wt \om(\la)$
$(t=\ov{t})$
exist and are finite almost everywhere. (One could use also theorem for analytic functions with positive real part
from \cite[p. 58]{Koos}.) In view of the representation \eqref{1'} of $\om$ and equalities \eqref{I1+}, \eqref{NC0} and \eqref{NC4},
it is easy to see that the function $\om_0(\la)=\om(\la)-\wt \om(\la)$ has the following property:
 \begin{align}  \label{E2-}  &
\lim_{\la \to \exp\{\I t\}} \Re \big(\om_0(\la)\big)=0 \quad {\mathrm{for}} \quad t\in (\ell_1,\, \ell_2).
 \end{align}
Here we used the relations \eqref{NC0-} and \eqref{NC3} as well. 
Hence, according to \eqref{E1} the following equality holds for the non-tangential limits:
 \begin{align}  \label{E2}  &
\lim_{\la \to \exp\{\I t\}} \Re \big(\om(\la)\big)=\lim_{\la \to \exp\{\I t\}} \Re\big(\wt \om(\la)\big) \quad {\mathrm{for}} \quad t\in (\ell_1,\, \ell_2).
 \end{align} 
 Therefore,  we  denote the limits in \eqref{E2} as $\Re\big(\om(\E^{\I t})\big)$ and $\Re \big(\wt \om(\E^{\I t})\big)$, and \eqref{E2} takes the form
 \begin{align}  \label{E2'}
\Re\big(\om(\E^{\I t})\big)=\Re \big(\wt \om(\E^{\I t})\big) \qquad ( t\in (\ell_1,\, \ell_2)).
 \end{align} 
Now, using the definition of $\wt \om$ in \eqref{E1} and Fatou's theorem (see, e.g., \cite[p. 39]{Nik})  we have $\Re \big(\wt \om(\E^{\I t})\big)=\mu^{\prime}(t)$
for $t\in (\ell_1,\, \ell_2)$.
Thus, almost everywhere on $(0,2\pi)$ equality  \eqref{E2'} yields 
\begin{align}  \label{E3}  &
\Re\big( \om(\E^{\I \th})\big)=\mu^{\prime}(\th).
 \end{align} 
 In view of \eqref{E3}, we rewrite \eqref{E0} in the form
\begin{align} \label{E4}
E(\vp,\wt \la)=-\frac{1}{4\pi}\int_{0}^{2\pi}|\E^{\I\th}-\wt \la |^{-2}(1-|\wt \la |^2)\ln \det\big(\Re\big( \om(\E^{\I \th})\big)\big)d \th.
\end{align}

Next, we turn to the linear fractional transformations \eqref{7}. According to \eqref{6},  \eqref{7}, and \eqref{7'} we have
\begin{align} \nn
\I\big(\vp(z)^*-\vp(z)\big)=& \big((c(z)P(z)+d(z)Q(z))^{-1}\big)^*
\\ \nn & \times 
\left(\begin{bmatrix} P(z)^* & Q(z)^*\end{bmatrix}
\begin{bmatrix} a(z)^* \\ b(z)^*\end{bmatrix}
\begin{bmatrix} c(z) & d(z)\end{bmatrix}
\begin{bmatrix} P(z) \\ Q(z)\end{bmatrix} \right.
\\ \nn &   \left.
+\begin{bmatrix} P(z)^* & Q(z)^*\end{bmatrix}
\begin{bmatrix} c(z)^* \\ d(z)^*\end{bmatrix}\begin{bmatrix} a(z) & b(z)\end{bmatrix}
\begin{bmatrix} P(z) \\ Q(z)\end{bmatrix}\right)
\\  \label{En}  & \times
(c(z)P(z)+d(z)Q(z))^{-1}
\\ \nn 
=& \big((c(z)P(z)+d(z)Q(z))^{-1}\big)^*
\\ \nn & \times
\left(\begin{bmatrix} P(z)^* & Q(z)^*\end{bmatrix}\clu(z)^*J\clu(z)
\begin{bmatrix} P(z) \\ Q(z)\end{bmatrix} \right)
\\ \nn & \times
(c(z)P(z)+d(z)Q(z))^{-1}.
\end{align} 
Recalling \eqref{En2-}, we see that
\begin{align} \label{En2}
\clu(\xi)^*J\clu(\xi)\equiv J \quad (\xi\in \BR).
\end{align}
Moreover, \eqref{NC} implies that
\begin{align} \label{Enn1}
z=(\ov{\up}\la - \up)\big/(\la-1),
\end{align}
and we set, correspondingly,
\begin{align} \label{Enn2}
\xi=\xi(\th)=(\ov{\up}\E^{\I \th} - \up)\big/(\E^{\I \th}-1).
\end{align}

Using \eqref{En}, \eqref{En2} and \eqref{Enn2}, we rewrite \eqref{E4} in the form
\begin{align} \nn
E(\vp,\wt \la)=&-\frac{1}{4\pi}\int_{0}^{2\pi}|\E^{\I\th}-\wt \la |^{-2}(1-| \wt \la |^2)
\\ & \label{E5} \qquad \quad
\times \Big(-2\ln\big|\det \big(c(\xi)P(\xi)+d(\xi)Q(\xi)\big)\big| 
\\ \nn & \qquad \qquad
\left. 
+\ln \det\left(\frac{1}{2}\begin{bmatrix} P(\xi)^* & Q(\xi)^*\end{bmatrix}J
\begin{bmatrix} P(\xi) \\ Q(\xi)\end{bmatrix}
\right)\right)d\th .
\end{align}

Let us study  first   the  case of the following nonsingular  pairs with the property-$J$: 
\begin{align} \label{Enn3}
P(z)=\psi(z)\in N_0=N, \quad Q\equiv \I I_p.
\end{align}
Each nonsingular  pair with the property-$J$, such that $\det Q(z)\not\equiv 0$, may be substituted by the
equivalent pair (generating the same $\vp(z)$) of the form \eqref{Enn3}. Using \eqref{Enn3}, we rewrite \eqref{E5} in the form
\begin{align} \label{E6}
E(\vp,\wt \la)=&
E(\psi, \wt \la)
\\ & \nn
+\frac{1}{2\pi}
\int_{0}^{2\pi}|\E^{\I\th}-\wt \la |^{-2}(1-|\wt \la |^2)\ln\big|\det \big(c(\xi)\psi(\xi)+\I d(\xi)\big)\big|d\th ,
\end{align}
where $\xi(\th)$ is defined in \eqref{Enn2}.
\begin{Rk}\label{RkAKr}
Formula \eqref{E6} is a generalization for the sign-indefinite case of the formula \cite[(11)]{ArK2}
$($see the corresponding Theorems 2 and 3 in \cite{ArK2}$)$.
\end{Rk}
The general case of the pairs $\{P,Q\}$ is dealt with in Remark \ref{RkContr}.
\begin{Rk}\label{RkCk} Our next considerations are similar to the corresponding considerations
in the proof of the  Szeg{\H o} limit theorem $($indefinite case$)$, see \cite[pp. 480-482]{ALSJFA}.
In particular, we use the fact that the entries of $\om(\la)-T(\la)$, where $\om \in C_{\vk}$ and $T$ is taken from
the representation \eqref{1'} of $\om$, belong to some Hardy class $H_{\de}$  $(\de>0)$.
This fact follows from the expressions \cite[(2.8)]{RoS2} for $\frac{\E^{\I\th}+\la}{\E^{\I\th}-\la}+L(a_i,\rho_i, \th,\la)$
and from the convergence inequality \eqref{NC0-}. As a result, for the entries $\om_{ij}$ of $\om$ we have for some
$\de >0$ the inequalities:
\begin{align} \label{Enn4}
\ov{\lim}_{r\to 1-0}\int_0^{2\pi}\big|\om_{ij}(r\E^{\I \th})\big|^{\de}d\th <\infty.
\end{align}
\end{Rk}
Further we assume that
\begin{align} \label{Enn5}
\clh=\BC^n,
\end{align}
and so $A$ and $S$ are $n\times n$ matrices.
\begin{Nn}\label{Nnq} Let us introduce the notation 
\begin{align} \label{Enn6}
q(\la)=  \left( \det\left(A^*-\frac{1}{z(\la)} I_n\right)\right)^p\det\big(c(z(\la))\psi(z(\la))+\I d(z(\la))\big)
\end{align}
where $z(\la)$ is given in \eqref{Enn1}. 

The number of different eigenvalues of $A$ in $\BC_+$
we denote by $\vt$.
\end{Nn}
\begin{Pn}\label{Pnzer} Let  $\clh=\BC^n$, $S\in \clp_{\vk}$, $\det S\not=0$ and \eqref{3} hold.
Assume that \eqref{9} is valid and $\ker \Phi_2=\{0\}$. Then, the number $\kappa$ of different zeros of $q(\la)$ in $\BD$ is finite.
Moreover, we have
\begin{align} \label{Enn7}
\kappa \leq \vk +\vt.
\end{align}
\end{Pn}
\begin{proof}. Similar to \cite{SaL0} and \cite[(1.17)]{ALSJFA} we introduce the matrix function
\begin{align} \label{Enn8}
A_{\psi}:=A-\Phi_2\big(\I \Phi_1^*-\psi\big(1\big/ \ov{z}\big)^*\Phi_2^*\big)S^{-1} \qquad (z\in \BC_+\cup \BC_-),
\end{align}
where 
\begin{align} \label{Enn8'}
\psi(z):= \psi(\ov{z})^* \quad {\mathrm{for}} \quad z\in \BC_-.
\end{align}
From \eqref{3} and \eqref{Enn8} we derive
\begin{align} \label{Enn9}
A_{\psi}S-SA_{\psi}^*=\Phi_2\big(\psi\big(1\big/ \ov{z}\big)^*-\psi\big(1\big/ \ov{z}\big)\big)\Phi_2^*.
\end{align}
In this proof we will consider $A_{\psi}$ and $\psi(z)$ for $z\in \BC_+$. Assume first that
\begin{align} \label{Enn10}
\psi \equiv {\mathrm{const}}, \quad \I(\psi^*-\psi)>0.
\end{align}
In this case, $A_{\psi}^*$ is $S$-dissipative  and does not have real eigenvalues. The first fact follows
from \eqref{Enn10} and the second fact is proved by contradiction. Indeed, assuming that
\begin{align} \label{Enn11}
A_{\psi}^*f=cf \quad (f\in \BC^n, \quad f\not=0, \quad c=\ov{c}),
\end{align}
and multiplying both parts of \eqref{Enn9} by $f^*$ from the left and by $f$ from the right, we obtain
\begin{align} \label{Enn12}
0=f^*\Phi_2\big(\psi^*-\psi\big)\Phi_2^*f
\end{align}
Recall that $\I(\psi^*-\psi)$ is sign-definite, and so equality \eqref{Enn12} yields $\Phi_2^*f=0$.
Hence, taking into account \eqref{Enn8} and \eqref{Enn11}, we derive
$A_{\psi}^*f=A^*f=cf$. Since (in view of \eqref{9}), $A^*$ does not have real eigenvalues,
we arrive to a contradiction.

We showed that $A_{\psi}^*$ is $S$-dissipative and does not have real eigenvalues. Recall also that $S\in \clp_{\vk}$. 
Thus, $A_{\psi}^*$ has $\vk$ eigenvalues in $\BC_-$ counting multiplicities. In the same way as in \cite[p. 477]{ALSJFA}
it follows that $\det(A_{\psi}-zI_n)$ has no more than $\vk$ zeros (counting multiplicities) in $\BC_+$ for all
$\psi(z) \in N_0$. We note that $\psi$ may now depend on $z$ and $\det(A_{\psi}-zI_n)$ means $\det(A_{\psi(z)}-zI_n)$.

According to \eqref{Enn8} we have the equality
\begin{align} \nn &
(A_{\psi}-z I_n)^{-1}\Phi_2-(A-z I_n)^{-1}\Phi_2 
\\ & \label{Enn13}
=(A_{\psi}-z I_n)^{-1}
\Phi_2\big(\I \Phi_1^*-\psi\big(1\big/ \ov{z}\big)^*\Phi_2^*\big)S^{-1} (A-z I_n)^{-1}\Phi_2,
\end{align}
which we rewrite in the form
\begin{align} & \label{Enn14}
(A_{\psi}-z I_n)^{-1}
\Phi_2 G(z)=
 (A-z I_n)^{-1}\Phi_2, 
\end{align}
where
\begin{align} \label{Enn15}
G(z)=I_p-\big(\I \Phi_1^*-\psi\big(1\big/ \ov{z}\big)^*\Phi_2^*\big)S^{-1} (A-z I_n)^{-1}\Phi_2.
\end{align}
On the other hand formula \eqref{6} implies that
\begin{align}\nn
c(z)\psi(z)+\I d(z)&=\I\big(I_p+\Phi_2^*(A^*-(1/z)I_n)^{-1}S^{-1}(\I \Phi_1+ \Phi_2\psi(z))\big)
\\ &  \label{Enn16}
=\I G\big(1\big/ \ov{z}\big)^*.
\end{align}
Finally, we rewrite \eqref{Enn14} in the form
\begin{align} & \label{Enn14'}
G\big(1\big/ \ov{z}\big)^*
\Phi_2^*\left(A_{\psi}^*-\frac{1}{z} I_n\right)^{-1}=
 \Phi_2^*\left(A^*-\frac{1}{z} I_n\right)^{-1}.
\end{align}

Recall  that the number
of different zeros of  $\det(A-zI_n)$ in $\BC_+$ equals $\vt$ and that $\det(A_{\psi}-zI_n)$ has no more than $\vk$ zeros (counting multiplicities) in $\BC_+$. 
Hence, the same holds for $\det\left(A^*-\frac{1}{z} I_n\right)$ and $\det\left(A_{\psi}^*-\frac{1}{z} I_n\right)$,
respectively. Taking into account relations \eqref{Enn15}, \eqref{Enn14'} and $\ker \Phi_2=\{0\}$, we see that  $G\big(1\big/ \ov{z}\big)^*$ is holomorphic and
does not have zeros in $\BC_+$ excluding, possibly zeros of $\det\left(A^*-\frac{1}{z} I_n\right)$ and $\det\left(A_{\psi}^*-\frac{1}{z} I_n\right)$.
Hence,  \eqref{Enn16} yields that the number of different zeros of 
$$\left( \det\left(A^*-\frac{1}{z} I_n\right)\right)^p\det\big(c(z)\psi(z)+\I d(z)\big)$$
in $\BC_+$ is less or equal $\vt+\vk$. Clearly, the same is valid, if we switch from $z$ to $z(\la)$, where $z(\la)$ is given by \eqref{Enn1},
and consider zeros of the function $q(\la)$ (see \eqref{Enn6}) in $\BD$. That is, \eqref{Enn7} is proved.
\end{proof}
\begin{Nn}\label{Nnzer} Different zeros of $q(\la)$ are denoted by $\la_1, \la_2, \ldots, \la_{\kappa}$
and their multiplicities are denoted by $\eta_1, \eta_2, \ldots, \eta_{\kappa}$, respectively.
\end{Nn}
It follows from \eqref{6} that $q(\la)$ is analytic in $\BD$. Moreover, since the entries of $\psi(z(\la))$ belong to $H_{\de}$, formula \eqref{6}
yields that $q(\la)$ belongs to some Hardy class $H_{\de}$ as well. Hence, using Proposition \ref{Pnzer} and Notation \ref{Nnzer} we
write down $q$ in the form
\begin{align}\label{Enn17}
q(\la)=B(\la)D(\la), \quad B(\la)=\prod_{i=1}^{\kappa}\big((\la -\la_i)\big/(1-\ov{\la_i}\la)\big),
\end{align}
where $B(\la)$ is the Blaschke product and $D(\la)$ belongs to $H_{\de}$ for some $\de>0$ and does not have zeros in $\BC_+$.

\begin{Tm}\label{TmEntr} Let \eqref{Enn3}   and \eqref{Enn5} hold, assume that the matrices  $A$, $S$, $\Phi_1$ and $\Phi_2$
satisfy the matrix identity \eqref{3}, that $A$ satisfies \eqref{9}, that $S\in \clp_{\vk_S}$, that $\det S\not=0$ and that $\ker \Phi_2=\{0\}$.

Then, the entropy functional on $\vp\in \cln(\clu)$ $($where $\clu$ is given by \eqref{6}$)$ satisfies the equality
\begin{align}\nn
E(\vp,\wt \la)=&E(\psi, \wt \la)
+\ln\big|\det\big(c(z(\wt \la))\psi(z(\wt \la))+\I d(z(\wt \la))\big)\big|
\\  \label{MF} &
+p\ln\left|  \det\left(A^*-\frac{1}{z(\wt \la)} I_n\right)\right|
-\ln|B(\wt \la)|
\\ & \nn
-\frac{p}{2\pi}
\int_{0}^{2\pi}|\E^{\I\th}-\wt \la |^{-2}(1-|\wt \la |^2)\ln\left|  \det\left(A^*-\frac{1}{\xi( \th)} I_n\right)\right|d\th 
\end{align}
$($if only the terms on the right-hand side above are finite$)$. Moreover, if we set $\psi(z):=\psi(\ov{z})^*$ for $z\in \BC_-$ and require additionally  that 
\begin{align}\label{adcond}
\det\big(c(z)\psi(z)+\I d(z)\big)\not=0 \quad {\mathrm{for}} \quad z \quad {\mathrm{such \,\,  that}}\quad 1/z\in\s(A),
\end{align}
 then the matrix functions $\vp$ are solutions of indefinite interpolations problems \eqref{23} and
 formula \eqref{MF} provides the values of entropy functionals for these solutions.
\end{Tm}
\begin{proof}.
Together with the matrix function  $c(z)\psi(z)+\I d(z)$ we consider the matrix function $\big(c(z)\psi(z)+\I d(z)\big)^{-1}$.
It is immediate from \eqref{6} and \eqref{7}  that
\begin{align} \label{E7}
\begin{bmatrix} \vp(z) \\ \I I_p \end{bmatrix}=\I \, \clu(z) \begin{bmatrix} \psi(z) \\ \I I_p \end{bmatrix}\big(c(z)\psi(z)+\I d(z)\big)^{-1}
\end{align}
Relation \eqref{En2-} implies that $\clu(z)^{-1}=J\clu(\ov{z})^*J$, and we rewrite \eqref{E7} in the form
\begin{align} \label{E8}
J\clu(\ov{z})^*J\begin{bmatrix} \vp(z) \\ \I I_p \end{bmatrix}= \begin{bmatrix} \I \psi(z) \\ - I_p \end{bmatrix}\big(c(z)\psi(z)+\I d(z)\big)^{-1}.
\end{align}
In particular, we have
\begin{align} \label{E9}
- \begin{bmatrix} I_p & 0 \end{bmatrix}\clu(\ov{z})^*\begin{bmatrix}  \I I_p \\ \vp(z) \end{bmatrix}= \big(c(z)\psi(z)+\I d(z)\big)^{-1}.
\end{align}

According to \eqref{NC}, \eqref{Enn6}, \eqref{Enn17} and \eqref{E9} we have
\begin{align} \label{Enn18}
D(\la)^{-1}= \frac{(-\I)^p  B(\la)}{\left( \det\left(A^*-\frac{1}{z(\la)} I_n\right)\right)^p} 
\det\left(\begin{bmatrix} I_p & 0 \end{bmatrix}\clu(\ov{z(\la)})^*\begin{bmatrix}   I_p \\ \om(\la) \end{bmatrix}\right).
\end{align}
Using \eqref{6}, \eqref{Enn4} and \eqref{Enn18} we obtain
\begin{align} \label{Enn19}
\ov{\lim}_{r\to 1-0}\int_0^{2\pi}\big|1\big/ D(r\E^{\I \th})\big|^{\de}d\th <\infty
\end{align}
for some $\de >0$. Formula \eqref{Enn19} means that $D(\la)^{-1}$ belongs 
some Hardy class, and we recall that $D(\la)$ belongs 
some Hardy class as well. Thus, $D(\la)$ is an outer function.
It is immediate from the parameter representation of the outer function (see, e.g. \cite[p. 76]{Koos})
that 
\begin{align} \label{Enn20}
\ln |D(\wt \la)|=\frac{1}{2\pi}
\int_{0}^{2\pi}|\E^{\I\th}-\wt \la |^{-2}(1-|\wt \la |^2)\ln\big|D(\E^{\I \th})\big|d\th .
\end{align}

Next, using  \eqref{Enn6} and \eqref{Enn17} we derive
\begin{align} \label{Enn21}
\det \big(c(\xi(\th))\psi(\xi(\th))+\I d(\xi(\th))\big)=\frac{B(\E^{\I \th})D(\E^{\I \th})}{ \left( \det\left(A^*-\frac{1}{\xi( \th)} I_n\right)\right)^p},
\end{align}
where $\xi(\th)$ is given in \eqref{Enn2}. We note that according
to \eqref{Enn17} the equality
\begin{align} \label{Enn22}
|B(\E^{\I \th}|=1
\end{align}
is valid. In view of \eqref{Enn20}--\eqref{Enn22}, we rewrite \eqref{E6} in the form
\begin{align} \label{E6'}
E(\vp,\wt \la)=&E(\psi, \wt \la)+\ln |q(\wt \la)|-\ln|B(\wt \la)|
\\ & \nn
-\frac{p}{2\pi}
\int_{0}^{2\pi}|\E^{\I\th}-\wt \la |^{-2}(1-|\wt \la |^2)\ln\left|  \det\left(A^*-\frac{1}{\xi( \th)} I_n\right)\right|d\th ,
\end{align}
Now, \eqref{Enn6} and \eqref{E6'} imply \eqref{MF}.

Taking into account \cite[Theorem 5.2]{RoS1} (for the case $p>1$), we see that the conditions
of Theorems \ref{TmInt} and \ref{TmIntp1} are fulfilled if \eqref{adcond} holds.  
In other words, if \eqref{adcond} holds, the matrix functions $\vp$ are solutions of indefinite interpolations problems \eqref{23}.
\end{proof}
\begin{Rk}\label{RkApsi} In the proof of Proposition \ref{Pnzer}, we studied $($starting from \eqref{Enn10}$)$ the case
of $z\in \BC_+$. We could substitute the inequality in \eqref{Enn10} with $\I(\psi-\psi^*)>0$, in which case $A_{\psi}^*$
would be $-S$-dissipative. Quite similar to the considerations after  \eqref{Enn10} one shows that $\det(A_{\psi(z)}-zI_n)$
has no more than $\vk$ zeros $($counting multiplicities$)$ in $\BC_-$. Here $\psi\in N_0$ $($in $\BC_+)$ and is introduced
via \eqref{Enn8'} in $\BC_-$. Clearly, relations \eqref{Enn14}--\eqref{Enn14'} hold in $\BC_+\cup \BC_-$.
\end{Rk}
\begin{Rk}\label{RkAll}
If one wants to consider entropy for all the solutions of indefinite interpolation problem
one may easily swithch from the matrix functions $\psi\in N$ to contractive matrix functions $\phi$.
See Remarm \ref{RkContr}, Definition \ref{DnContr} and formula \eqref{E5'} as well as the considerations
after \eqref{E5'}.
\end{Rk}
\section{Entropy functionals and indefinite \\ Caratheodory problem}\label{Car}
\setcounter{equation}{0}
Consider the case of block Toeplitz matrices, that is, the case
\begin{align}\label{c0}
S=S(n)=\{s_{j-i}\}_{i,j=1}^n \in \clp_{\vk}, \quad \det S\not=0,
\end{align}
where $s_k$ are $p\times p$ blocks. Block Toeplitz matrices $S$ are unique solutions
of the matrix identities (see \cite{ALS73, SaA80, ALSJFA}):
\begin{align}& \label{c1}
AS-SA^*= \I \Pi J \Pi^*;  \quad \Pi=\begin{bmatrix}\Phi_1 & \Phi_2\end{bmatrix},  
\end{align}
where
\begin{align} & \label{c2}
A= \left\{ a_{j-i}^{\,} \right\}_{i,j=1}^n,
          \quad  a_k  =  \left\{ \begin{array}{lll}
                                  0 \, & \mbox{ for }&
k> 0  \\
\displaystyle{\frac{\I}{ {\, 2 \,}}} \,
                                I_p
                                  & \mbox{ for }& k =
0  \\
                                \, \I \, I_p
                                  & \mbox{ for }& k <
0
                          \end{array} \right.
, \qquad J=\begin{bmatrix} 0 & I_p \\ I_p & 0 \end{bmatrix};
\\ & \label{c3}
\Phi_2 = \left[
\begin{array}{c}
I_{p}  \\ I_p \\ \cdots \\ I_{p}
\end{array}
\right], \quad \Phi_1 =  \left[
\begin{array}{l}
s_0/2  \\ s_0/2 + s_{-1} \\ \cdots \\ s_0/2+ s_{-1} + \ldots +
s_{1-n}
\end{array}
\right] +\I \Phi_2 \nu, \quad \nu = \nu^*;
\\ & \label{c3'}
A=A(n), \quad \Pi=\Pi(n), \quad \Phi_1=\Phi_1(n), \quad \Phi_2=\Phi_2(n).
\end{align}
 For related results and discussions on the operator identities or operators and matrices with the displacement structure see, for instance,
 \cite{KaS, KaS2, SaSaR,
SaL1, SaL1'}.  The next proposition is immediate from \eqref{c0}--\eqref{c3}.
\begin{Pn} \label{PnTriple} The matrices $A$, $S$, $\Phi_1$ and $\Phi_2$ given by \eqref{c0}, \eqref{c2} and \eqref{c3} satisfy
conditions on these matrices from Theorems \ref{TmInt} and \ref{TmEntr}, and the corresponding space $\clh$, where $A$ and 
$S$ act, is finite dimensional. More precisely,
$\clh=\BC^{pn}$.
\end{Pn}
Let us {\it add index $n$ in the notation $\clp_{\vk}$} and write $S(n)\in \clp_{\vk,n}$ instead of $S(n)\in \clp_{\vk}$
when we consider $pn\times pn$ matrices. We say that $\{s_{j-i}\}_{i,j=1}^{\infty} \in \clp_{\vk,\infty}$ if all
the reductions $\{s_{j-i}\}_{i,j=1}^{\wh n}$ $(n_0 \leq \wh n <\infty)$ starting from some $n_0$ belong to $\clp_{\vk, \wh n}$.
\begin{Dn}\label{DnExt} Matrices $S(\wt n)=\{s_{j-i}\}_{i,j=1}^{\wt n}$ $(n<\wt n\leq \infty)$ 
are called extensions of the matrix   $S=S(n)=\{s_{j-i}\}_{i,j=1}^n$.
\end{Dn}
The following result  may be derived as a special case of Theorem \ref{TmInt}
and presents a reformulation of \cite[Theorem 4.1]{ALSJFA}.
\begin{Tm}\label{TmExt}
Let the matrices  $A$,  $\Phi_1$ and $\Phi_2$ by given by \eqref{c2} and \eqref{c3},
let $S$ be given by the equalities in \eqref{c0}, and let the relations $S\in \clp_{\vk,n}$
and $\det S\not=0$ hold. Assume that the matrix functions $\vp(z,P,Q)$ are given
by \eqref{7}, where $\clu$ is given by \eqref{6} and the pairs $\{P,\,Q\}$ are nonsingular pairs  with the property-$J$
satisfying inequality
\begin{align} \label{c4}
\det\left(\begin{bmatrix} 0 & \ldots & 0 & I_p \end{bmatrix}S^{-1}\Pi J \begin{bmatrix} P(2\I) \\ Q(2\I)  \end{bmatrix}\right)\not=0.
\end{align}
Then the functions $-\I \vp\left(2\I \frac{1-\la}{1+\la}\right)$ admit Taylor expansions with the same first $n$ Taylor coefficients
determined by $S=\{s_{j-i}\}_{i,j=1}^n:$
\begin{align} \label{c5}
-\I\vp\left(2\I \frac{1-\la}{1+\la}\right)=\left(\frac{s_0}{2}-\I \nu\right)+s_1 \la + \ldots +s_{n-1}\la^{n-1}+\ldots .
\end{align}
Moreover, further Taylor coefficients of the functions $\vp$ generate extensions $S(\wt n)=\{s_{j-i}\}_{i,j=1}^{\wt n}$ 
of $S$ belonging to the classes $\clp_{\vk,\wt n}:$
\begin{align} \nn
-\I\vp\left(2\I \frac{1-\la}{1+\la}\right)=&\left(\frac{s_0}{2}-\I \nu\right)+s_1 \la + \ldots +s_{n-1}\la^{n-1}+s_{n}\la^{n} +\ldots
\\ \label{c6} &
+s_{\wh n}\la^{\wh n}+\ldots 
\end{align}
$(n<\wh n \leq \wt n)$, and all the extensions of $S$ belonging to the classes $\clp_{\vk,\wt n}$ are generated in this way.
\end{Tm}
While reformulating \cite[Theorem 4.1]{ALSJFA} as Theorem \ref{TmExt} we took into account that the coefficients
$\{w_{ij}(\la)\}_{i,j=1}^n$ (see \cite[(1.6)]{ALSJFA}) of the linear fractional transformation $\vp_{\a}(\la)$ of the form
\cite[(1.8)]{ALSJFA} and $\clu(z)$  given by \eqref{6} are connected by the relation:
\begin{align} \label{c7}
J\left(\left\{w_{ij}\left(\frac{1-(\I/2)\ov{z}}{1+(\I/2)\ov{z}}\right)\right\}_{i,j=1}^n \right)^*J=\clu(z),
\end{align}
Formula \eqref{c7} implies the equalities $\vp_{\a}(\ov{\la})^*=-\I \vp\left(2\I \frac{1-\la}{1+\la},P,Q\right)$,
where the matrix functions $\vp_{\a}(\la)$ are given by \cite[(1.8)]{ALSJFA}, the matrix functions $\vp(z,P,Q)$ are given
by \eqref{7}, and there is a simple one to one correspondence between the matrix functions $\a(\la)$
and nonsingular pairs $\{P(z), \, Q(z)\}$ with the property-$J$.

Next, let us consider condition \eqref{c4} in greater detail. We start with a simple auxiliary lemma.
\begin{La}\label{LaAux} Let $f_1, f_2$ and $g_1, g_2$ be $p \times 2p$ matrices, 
let all four matrices have rank $p$ and assume that the equalities
\begin{align} \label{A1}
f_1g_1^*=f_2g_2^*=0
\end{align}
hold. Then, the inequalities
\begin{align} \label{A2}
\det (f_1f_2^*)\not=0 \quad {\mathrm{and}}\quad  \det (g_1g_2^*)\not=0
\end{align}
are equivalent. 
\end{La}
\begin{proof}.  Let us show that the first inequality in \eqref{A2} yields the second (clearly, the fact that the second
inequality yields the first is proved in the same way). Indeed, if $\det (g_1g_2^*)=0$ there is $h\in \BC^p, \,\, h\not=0$
such that $hg_1g_2^*=0$. Hence, in view of $f_2g_2^*=0$, there is $\wh h\not=0$ such that $hg_1=\wh h f_2$,
and therefore $f_1g_1^*=0$ implies that $f_1f_2^*\wh h^*=0$, which contradicts $\det (f_1f_2^*)\not=0$.
\end{proof}

\begin{La}\label{LaDet} Let the matrices  $A$,  $\Phi_1$ and $\Phi_2$ by given by \eqref{c2} and \eqref{c3},
let $S$ be given by the equalities in \eqref{c0}, and let the relations $S\in \clp_{\vk,n}$
and $\det S\not=0$ hold. Assume that the pair $\{P,\, Q\}$ is given
by \eqref{Enn3}. Then, condition \eqref{c4} is equivalent to the condition
\begin{align} \label{c9}
\det\big(c(-2\I)\psi(2\I)^*+\I d(-2\I)\big)\not=0.
\end{align}
\end{La}
\begin{proof}. According to \cite[p. 452]{ALSJFA} we have
\begin{align} \label{c10}
\left(A-\frac{1}{z}I_{pn}\right)^{-1}\Phi_2=-\frac{z}{1-(\I/2) z}\begin{bmatrix} I_p \\ \frac{1+(\I/2) z}{1-(\I/2) z} I_p \\ \ldots 
\\ \left(\frac{1+(\I/2) z}{1-(\I/2) z}\right)^{n-1} I_p \end{bmatrix}.
\end{align}
Setting
\begin{align} \label{c8}
Y=Y(n):=\begin{bmatrix} 0 & \ldots & 0 & I_p \end{bmatrix}S(n)^{-1}\Pi (n)
\end{align}
and using \eqref{6} and \eqref{c10}, we see that
\begin{align} \label{c11}
\lim_{z\to -2\I}\big(1-(\I/2) z\big)^{n}\begin{bmatrix} c(\ov{z}) & d(\ov{z})  \end{bmatrix}^*=2^n JY^*.
\end{align}
Moreover, formulas \eqref{6}, \eqref{En2-} and \eqref{c11} imply the equality
\begin{align} \label{c12}
\begin{bmatrix} c({-2\I}) & d({-2\I})  \end{bmatrix}Y^*=0,
\end{align}
where, in view of  \eqref{6} and \eqref{c10}, we have
\begin{align} \label{c13}
\begin{bmatrix} c({-2\I}) & d({-2\I})  \end{bmatrix}=\begin{bmatrix} 0 & I_p  \end{bmatrix} - \begin{bmatrix}I_p & 0 &\ldots & 0  \end{bmatrix}S^{-1}\Pi J.
\end{align}
We will need the equalities
\begin{align} \label{c14}
{\mathrm{rank}} \begin{bmatrix} c({-2\I}) & d({-2\I})  \end{bmatrix}={\mathrm{rank}} \, Y=p.
\end{align}
In order to show that rank $Y=p$, we partition the following matrices into the $p\times p$ blocks:
\begin{align} \label{c15}&
\begin{bmatrix} 0 & \ldots & 0 & I_p \end{bmatrix}S^{-1}=\begin{bmatrix} t_{n1} & t_{n2}& \ldots & t_{nn} \end{bmatrix},
\\ & \label{c16}
 YJ\Pi^*S^{-1}=\begin{bmatrix} q_{n,1} & q_{n,2}& \ldots & q_{n,n} \end{bmatrix}.
\end{align}
It easily follows from \eqref{c1} and is immediate from \cite[(12)]{SaA80} that
\begin{align} \label{c17}
t_{nk}=q_{n,k}-q_{n,k+1} \quad (q_{n,n+1}:=0).
\end{align}
Thus, if rank $Y\not=p$, we have $hY=0$ for some $h\in \BC^p, \,\, h\not=0$,
and formulas \eqref{c15}--\eqref{c17} yield
\begin{align} \label{c18}&
\begin{bmatrix} 0 & \ldots & 0 & h \end{bmatrix}S^{-1}=0 ,\end{align}
which contradicts the invertibility of $S$. Hence, indeed, rank $Y=p$ .

Using the equality rank $Y=p$, one can show that there is an extension $S(n+1) \in \clp_{\vk, n+1}$ such that
$\det S(n+1)\not=0$ (see, e.g., \cite[Lemma 8]{SaA80}). Now, relations \eqref{c13} and \cite[(2.16)]{ALSJFA}
yield the equality 
$${\mathrm{rank}}\begin{bmatrix} c({-2\I}) & d({-2\I})  \end{bmatrix}=p,$$ 
and \eqref{c14} is proved.

Finally we set $$f_1=Y, \,\, f_2= \begin{bmatrix} -\I I_p & \psi({2\I})^* \end{bmatrix},\,\, g_1=\begin{bmatrix} c({-2\I}) & d({-2\I})  \end{bmatrix}, \,\,
g_2= \begin{bmatrix} \psi({2\I}) & -\I I_p  \end{bmatrix}.$$
Taking into account \eqref{c12} and \eqref{c14}, we see that $f_1,f_2,g_1$ and $g_2$ satisfy conditions
of Lemma \ref{LaAux}. It follows that the inequalities 
\begin{align} \label{c4'}
\det\left(\begin{bmatrix} 0 & \ldots & 0 & I_p \end{bmatrix}S^{-1}\Pi J \begin{bmatrix} \psi(2\I) \\ \I I_p  \end{bmatrix}\right)\not=0
\end{align}
and \eqref{c9} are equivalent.
\end{proof}
In order to derive the entropy formula and include the important case $\wt \la =0$ one could remove
singularities of $\det\big(c(z(\la))\psi(z(\la))+\I d(z(\la))\big)$ inside $\BD$ in a simpler way than in \eqref{Enn6} and 
introduce
$q(\la)$ as the product:
\begin{align} \label{Enn6'}
q(\la)=   \det\left(A^*-\frac{1}{z(\la)} I_{pn}\right)\det\big(c(z(\la))\psi(z(\la))+\I d(z(\la))\big).
\end{align}
We note that in the case \eqref{c0}--\eqref{c3} of Toeplitz matrices (and in view of \eqref{Enn1},  where $\up=2\I$)
the following 
equality holds:
\begin{align}& \label{c21}
\det\left(A^*-\frac{1}{z( \la)}I_{pn}\right)=\left(\frac{ \la}{\I( \la +1)}\right)^{pn}.
\end{align}
Hence, it will be even more convenient to consider $\wt q(\la)$ of the form
\begin{align} \label{Enn6!}
\wt q(\la)=  \det\Big(\la^n\big(c(z(\la))\psi(z(\la))+\I d(z(\la))\big)\Big)
\end{align}
instead of $q(\la)$ given by \eqref{Enn6'}.
Now,  from \eqref{6}, \eqref{Enn1}, \eqref{c10} and \eqref{Enn6!} we derive
\begin{align}  \label{c25}
\wt q(\la)=&\det\left(\Big(\begin{bmatrix}0 &\la^n I_p  \end{bmatrix}\right.
\\ & \nn \left.
-({\la+1})
\begin{bmatrix}\la^{n-1} I_p & -\la^{n-2}I_p & \ldots &\left(-1\right)^{n-1}I_p   \end{bmatrix}S^{-1}\Pi J\Big)\begin{bmatrix}\psi(z(\la)) \\ \I I_p  \end{bmatrix}\right).
\end{align}
We introduce $\wt B(\la)$ as the Blashke product in the representation of  $\,\wt q(\la)$,
that is, $\wt B(\la)$ in terms of zeros $\la_i\in \BD$ of $\wt q(\la)$ $($counting multiplicities$)$ is given by the formula
\begin{align} \label{c29}
 \wt B(\la)=\prod_{i=1}^{\wt \kappa}\big((\la -\la_i)\big/(1-\ov{\la_i}\la)\big).
\end{align}
In a similar to the proof of Theorem \ref{TmEntr} way one can show that 
\begin{align} \label{c32}
\wt q(\la)=\wt B(\la)D(\la),
\end{align}
where $D(\la)$ is an outer function.
\begin{Tm}\label{ToepEntr} Let the matrices  $A$,  $\Phi_1$ and $\Phi_2$ be given by \eqref{c2} and \eqref{c3},
let $S$ be given by the equalities in \eqref{c0}, and let the relations $S\in \clp_{\vk,n}$
and $\det S\not=0$ hold. Assume that the pair $\{P, \, Q\}$ has the form $\{\psi, \I I_p\}$ $($as in \eqref{Enn3}$)$ and that \eqref{c4'} is valid.
Set $\up=2\I$ in \eqref{NC}, \eqref{Enn1} and \eqref{Enn2}.
Then, $\vp(z, \psi, \I I_p)$ defined by \eqref{7} is a solution of the interpolation problem \eqref{23} $($as well as the indefinite Caratheodory problem
described in Theorem \ref{TmExt}$)$
and its  entropy   is given by the formula:
\begin{align} \label{c28} &
E(\vp,\wt \la)=E(\psi, \wt \la)
+\ln\big|\wt q(\wt \la)\big|-\ln|\wt B(\wt \la)| ,
\end{align}
where $\wt q$ and $\wt B$ are given by \eqref{Enn6!} and \eqref{c29}, respectively.
\end{Tm}
\begin{proof}. It is easy to see that  the conditions of Theorem \ref{TmExt} are fulfilled. According to Lemma \ref{LaDet}, inequality \eqref{c4'}
is equivalent to \eqref{c9}, which in our case coincides with \eqref{adcond}. Thus, $\vp$ satisfies the conditions of Theorem \ref{TmEntr}
as well. Therefore, $\vp$ is a solution of the interpolation problem \eqref{23} and  of the indefinite Caratheodory problem
described in Theorem \ref{TmExt}. 

Taking into account that $\big|\wt B\big(\E^{\I \th}\big)\big|=1$, that $|\la |=1$ and that equalities \eqref{Enn6!} and \eqref{c32} hold, in the same way as in the proof of \eqref{MF}, we rewrite \eqref{E6} in the form \eqref{c28}.
\end{proof}
Formulas \eqref{c8} and  \eqref{c25} yield
\begin{align} \label{c45} &
\wt q(0)=(-1)^{p+n-1}\det\left(Y\begin{bmatrix}\I I_p  \\  \psi(2\I)   \end{bmatrix}\right).
\end{align}
\begin{Rk}\label{Rk0} Relations \eqref{c4'} and \eqref{c45} imply that $\wt q(\la)$ is well-defined  at $\la =0$ and moreover
$\wt q(0)\not=0$. Since $\wt q(0)\not=0$ and $\wt B(0)\not=0$, one can consider entropy $E(\vp, \wt \la)$ at $\wt \la=0$.
\end{Rk}

\begin{Rk}\label{RkContr} Sometimes, it is more convenient to use  contractive $p \times p$  matrix functions $\phi$
instead of the nonsingular pairs $\{P,Q\}$ with the property-$J$. For this purpose, one introduces $2p\times 2p$ matrix 
$\clw$ with the following property:
\begin{align} \label{c35} &
\clw:=\frac{1}{\sqrt{2}}\begin{bmatrix}I_p & -I_p \\  I_p & I_p   \end{bmatrix}, \quad \clw^* J \clw=j, \quad j:=\begin{bmatrix}I_p & 0 \\ 0   & -I_p   \end{bmatrix}.
\end{align}
We set:
\begin{align} \label{c36} &
\begin{bmatrix}\wh P(z) \\ \wh Q(z)  \end{bmatrix}:=\clw^{-1} \begin{bmatrix} P(z) \\  Q(z)  \end{bmatrix}, \quad \phi(z):=\wh Q(z) \wh P(z)^{-1}.
\end{align}
From \eqref{c35} and \eqref{c36}, it is easy to see that 
$$ {\mathrm{rank}} \begin{bmatrix}\wh P(z) \\ \wh Q(z)  \end{bmatrix}=p, \quad
 \begin{bmatrix}\wh P(z)^* & \wh Q(z)^*  \end{bmatrix}j \begin{bmatrix}\wh P(z) \\ \wh Q(z)  \end{bmatrix} \geq 0,$$
and so
$\wh P(z)$ is invertible and $\phi(z)$ is contractive when the  pair $\{P(z),Q(z)\}$ is nonsingular with the property-$J$.
Moreover, relations \eqref{7}, \eqref{c35} and \eqref{c36} imply that
\begin{align} \label{c37} &
\vp(z,P,Q)=\wh \vp(z,\phi),
\end{align}
where
\begin{align} \label{c38} &
\wh \vp(z,\phi)=\I(\wh a(z)+\wh b(z)\phi(z))(\wh c(z)+\wh d(z)\phi(z))^{-1}, \\
 \label{c39} &
 \wh \clu(z)= \begin{bmatrix}\wh a(z) & \wh b(z) \\ \wh  c(z) & \wh d(z) \end{bmatrix}:=
\clu(z)\clw .
\end{align}
Clearly, inequality  \eqref{7'} is equivalent to the inequality 
\begin{align} \label{c40} &
\det(\wh c(z)+\wh d(z)\phi(z))\not\equiv 0, 
 \end{align}
 and so $\cln(\clu)$ $($see Notation \ref{LFT}$)$ coincides with the set of functions $\wh \vp(z,\phi)$,
 where $\phi$ are contractive in $\BC_+$ and \eqref{c40} holds.
\end{Rk}
\begin{Dn}\label{DnContr} According to \cite[(1)]{ArK2}, the entropy $\wh E$ of the $p\times p$ matrix function
$g(\la)$ $($which is contractive in $\BD)$ is given by the formula
\begin{align} \label{c41} &
\wh E(g,\wt \la)=-\frac{1}{4\pi}\int_0^{2\pi}|\E^{\I\th}-\wt \la |^{-2}(1-| \wt\la |^2)\ln \det\big(I_p-g(\E^{\I\th})^*g(\E^{\I\th})\big)d\th .
 \end{align}
\end{Dn}
Taking into account \eqref{c35}, \eqref{c36}, \eqref{c39} and \eqref{c41} we rewrite
\eqref{E5} in the form
\begin{align}  \label{E5'} 
E(\vp,\wt \la)=&\wh E\big(\phi(z(\la)),\wt \la \big)+(p\ln 2)/2
\\ \nn
&+\frac{1}{2\pi}\int_{0}^{2\pi}|\E^{\I\th}-\wt \la |^{-2}(1-| \wt \la |^2)\ln\big|\det \big(\wh c(\xi)+\wh d(\xi)\phi(\xi)\big)\big| d\th ,
\end{align}
where $\xi(\th)$ is given by \eqref{Enn2}. 

In the case of indefinite  Caratheodory problem, we (similar to \eqref{Enn6!} and \eqref{c29}) set
\begin{align}  \label{c42} &
\wh q(\la)=  \det\Big(\la^n\big(\wh c(z(\la))+ \wh d(z(\la))\phi(z(\la))\big)\Big), \\
 \label{c43} &
 \wh B(\la)=\prod_{i=1}^{\wh \kappa}\big((\la -\la_i)\big/(1-\ov{\la_i}\la)\big),
\end{align}
where $\la_i\in \BD$ are zeros of $\wh q(\la)$ $($counting multiplicities$)$. Next, we factorise
$\wh q$: $\wh q(\la)=\wh B(\la)D(\la)$ and
in the same way as \eqref{Enn18} we obtain the equality
\begin{align} \label{Enn18'}
D(\la)^{-1}= \frac{\wh B(\la)}{\la^{pn}} 
\det\left(\begin{bmatrix} I_p & 0 \end{bmatrix}\wh \clu(\ov{z(\la)})^*\begin{bmatrix}   I_p \\ \om(\la) \end{bmatrix}\right).
\end{align}
Using \eqref{Enn4} and \eqref{Enn18'}, one can see that $D(\la)$ is again an outer function.
Hence, taking into account \eqref{E5'}, we have a somewhat stronger version of Theorem \ref{ToepEntr}.
\begin{Tm}\label{ToepEntr'} Let the matrices  $A$,  $\Phi_1$ and $\Phi_2$ be given by \eqref{c2} and \eqref{c3},
let $S$ be given by the equalities in \eqref{c0}, and let the relations $S\in \clp_{\vk,n}$
and $\det S\not=0$ hold. Set $\up=2\I$ in \eqref{NC}, \eqref{Enn1} and \eqref{Enn2}. Assume that the $p \times p$ matrix functions
$\phi(z)$ are contractive in $\BC_+$ and satisfy the inequality
\begin{align} \label{c4''}
\det\left(\begin{bmatrix} 0 & \ldots & 0 & I_p \end{bmatrix}S^{-1}\Pi J \clw \begin{bmatrix} I_p \\ \phi(2\I)  \end{bmatrix}\right)\not=0.
\end{align}
Then, the set of matrix functions $\wh \vp(z,\phi)$ of the form  \eqref{c38} describes $($via Taylor coefficients of $\wh \vp(z(-\la),\phi)$ as in Theorem \ref{TmExt}$)$ all the  solutions of  the indefinite Caratheodory problem.
The entropy functional on these matrix functions $\wh \vp$  is given by the formula:
\begin{align}  \label{E5''} 
E(\wh \vp,\wt \la)=&\wh E\big(\phi(z(\la)),\wt \la \big)+(p\ln 2)/2+\ln \big|\wh q(\wt \la)\big|-\ln \big|\wh B(\wt \la)\big|.
\end{align}
\end{Tm}
\section{Nonclassical Szeg{\H{o}} limit formula}  \label{lim}
\setcounter{equation}{0}
In order to introduce the entropy functionals, we transformed  functions $\vp(z)\in N_{\vk}$
into the functions $\om(\la)$ belonging to the class $C_{\vk}$. It is easy to see that the relation 
$\om(\la)\in C_{\vk}$ yields $\om(-\la)\in C_{\vk}$. In precisely the same way as before, one can
deal with the entropy $E_{\star}$ generated by the transformation of $\vp$ into  $\om_{\star}(\la):=\om(-\la)$.
In particular, relation \eqref{E4} takes the form
\begin{align} \label{S1}
E_{\star}(\vp,\wt \la)=-\frac{1}{4\pi}\int_{0}^{2\pi}|\E^{\I\th}-\wt \la |^{-2}(1-|\wt \la |^2)\ln \det\big(\Re\big(\om_{\star}(\E^{\I \th})\big)\big)d \th.
\end{align}
It is convenient to consider $E_{\star}$ instead of $E$ in the case of  Caratheodory problem. Indeed setting $\up=2\I$ in \eqref{Enn1}
we obtain
\begin{align} \label{S2}
\om_{\star}(\la):=\om(-\la)=-\I\vp\left(2\I \frac{1-\la}{1+\la}\right),
\end{align}
that is $\om_{\star}$ coincides with the function on the left-hand sides of important formulas \eqref{c5} and \eqref{c6}.
Using \eqref{S1} and the same arguments as in the proof of Theorem \ref{ToepEntr}, one rewrites \eqref{c28} in the form
\begin{align} \label{c28'} &
E_{\star}(\vp,\wt \la)=E_{\star}(\psi, \wt \la)
+\ln\big|\wt q(-\wt \la)\big|-\ln|\wt B(-\wt \la)| .
\end{align}
Now, let us fix $n$, $S(n)\in\clp_{\vk, n}$ and $\psi(z)\in N_0$. In the same way as $\vp$ generates block Toeplitz matrices $S(i)$ $(i>0)$
using Taylor coefficients from \eqref{c6}, the matrix function $\psi$ generates block Toeplitz matrices $\breve S(i)=\{\breve s_{k-j}\}_{j,k=1}^i$ $(i>0)$.
Since $\psi \in N_0$, we have $\breve S(i)\geq 0$ for all $i>0$. The notations, which we  introduce using $S(i)$ (e.g., $\Pi(i)$ and $Y(i)$) will obtain
an accent "breve" if we substitute $S(i)$ with $\breve S(i)$ (and we will write $\breve \Pi(i)$ and $\breve Y(i)$ in that case).
Introduce also the notations:
\begin{align} \label{c46} &
\det S(i)=\Lam_i, \qquad  \det \breve S(i)=\breve \Lam_i .
\end{align}
When $\breve S(i)> 0$,  the famous first Szeg{\H{o}} limit formula is valid:
\begin{align}\nn
\lim_{i\to \infty} \frac{\breve \Lam_i}{\breve \Lam_{i-1}}
&=\exp\left\{\frac{1}{2\pi}\int_0^{2\pi}\ln \big(\det(f(\th))\big)d\th\right\}, 
\\ \label{c47} &
=\exp\left\{\frac{1}{2\pi}\int_0^{2\pi}\ln \big(\det(f(2\pi - \th))\big)d\th\right\} 
\\
\label{c48} &
 f(\th):=\I\left(\psi\left(2\I\frac{\E^{\I \th}-1}{\E^{\I \th}+1}\right)^*
-\psi\left(2\I\frac{\E^{\I \th}-1}{\E^{\I \th}+1}\right)\right).
\end{align}
On the other hand, formulas \eqref{S1} and \eqref{S2} imply that
\begin{align} \label{S3}
E_{\star}(\psi,0)=-\frac{1}{4\pi}\int_{0}^{2\pi}\ln \det\big(f(2\pi -\th)/2\big)\big)d \th.
\end{align}
According to \eqref{c47} and \eqref{S3} we have
\begin{align}\label{S4}&
\lim_{i\to \infty} \frac{\breve \Lam_i}{\breve \Lam_{i-1}}
=2^p\exp\left\{-2E_{\star}(\psi,0) \right\}.
\end{align}
The following important correspondences exist between the matrices generated by $S(n+i)$ and $\breve S(i)$ (see \cite[(3.5), (3.6), (3.40)]{ALSJFA})
\begin{align} \label{S5}
\breve Y(i)^*=Y(n+i)^*C; \qquad \breve t_{ii}=C^*t_{n+i,n+i}C ,
\end{align}
where the matrix $C$ does not depend on $i$ and $\breve t_{ii}$ ($t_{n+i,n+i}$) is the right lower block of $\breve S(i)^{-1}$ (of $S(n+i)^{-1}$):
$$\breve t_{ii}:=\begin{bmatrix}0 & \ldots & 0 &I_p\end{bmatrix}\breve S(i)^{-1}\begin{bmatrix}0 & \ldots & 0 &I_p\end{bmatrix}^*.$$
From \eqref{c45} and \cite[(4.39)]{ALSJFA} it follows that 
\begin{align} \label{S6}
\big| \wt q(0)\big|^2=1\big/\det(C^*C).
\end{align}
Moreover, we have (see, e.g., \cite[(4.31)]{ALSJFA}):
\begin{align} \label{S7}
\det (\breve t_{ii})=   \frac{\breve \Lam_{i-1}}{\breve \Lam_i }    , \quad \det(t_{n+i,n+i})= \frac{ \Lam_{n+i-1}}{ \Lam_{n+i} }.
\end{align}
Relations \eqref{S4}--\eqref{S7} yield
\begin{align}\label{S8}&
\lim_{i\to \infty} \frac{ \Lam_i}{ \Lam_{i-1}}
=2^p\exp\left\{-2E_{\star}(\psi,0) \right\}\big| \wt q(0)\big|^{-2}.
\end{align}
Using the entropy formula \eqref{c28'} (at $\wt \la =0$) and taking into account \eqref{c29}, we rewrite
\eqref{S8} in the form
\begin{align}\label{S9}&
\lim_{i\to \infty} \frac{ \Lam_i}{ \Lam_{i-1}}
=2^p\exp\left\{-2E_{\star}(\vp,0) \right\}\big| \wt B(0)\big|^{-2}=2^p\exp\left\{-2E_{\star}(\phi,0) \right\}\prod_{j=1}^{\wt \kappa}|\la_j|^{-2}.
\end{align}
From the formula above, in view of \eqref{S1} and \eqref{S2} we derive  the nonclassical  Szeg{\H{o}} limit formula for matrices
$S(n)\in \clp(n)$ generated by $\vp\left(2\I \frac{1-\la}{1+\la}\right)$:
\begin{align}\nn &
\lim_{i\to \infty} \frac{ \Lam_i}{ \Lam_{i-1}}
=2^p\prod_{j=1}^{\wt \kappa}|\la_j|^{-2}\exp\left\{\frac{1}{2\pi}\int_0^{2\pi}\ln \left(\det \Im\left(\vp\left(2\I \frac{1-\E^{\I \th}}{1+\E^{\I \th}}\right)\right)\right)d\th\right\},
\end{align}
where $\la_j$ are the poles of $\det \vp\left(2\I \frac{1-\la}{1+\la}\right)$ counting multiplicities. For further details
see \cite[Section 4]{ALSJFA}.
On the asymptotics of determinants in some other important nonclassical cases see, for instance, \cite{BasE, Ehr}.

\bigskip

\noindent{\bf Acknowledgments.}
 {This research    was supported by the
Austrian Science Fund (FWF) under Grant  No. P29177.}


\begin{flushright}

I. Roitberg, \\
 e-mail: {\tt  	innaroitberg@gmail.com}

\vspace{0.5em} 

A.L. Sakhnovich,\\
Faculty of Mathematics,
University
of
Vienna, \\
Oskar-Morgenstern-Platz 1, A-1090 Vienna,
Austria, \\
e-mail: {\tt oleksandr.sakhnovych@univie.ac.at}

\end{flushright}

\end{document}